\documentstyle[12pt]{article}

\title{Multilinear quantum Lie operations}

\author{V.K. Kharchenko$^{\ast}$}  

\date{}  

\begin{document}   

\maketitle

\renewcommand{\thefootnote}{\fnsymbol{footnote}}

\setcounter{footnote}{-1}

\footnote{$\hspace*{-6mm}^{\ast}$

{\bf Pablished in Russian in Zap. Nauchn. Sem. S.-Petersburg. Otdel. Mat. Inst. Steklov. (POMI),
272 (2000), Vopr. Teor. Predst. Algebr. i Grupp. 7, 321-340, 351.}\\ 

Supported by 
CONACyT, M\'exico, grant 32130-E, and  PAPIIT UNAM, grant IN 102599.}

\renewcommand{\thefootnote}{\arabic{footnote}}\

\newcounter{par}

\newcounter{nom}[par]

\newtheorem{nom}[par]{Nom}

\newtheorem{theorem}{\bf Theorem}

\newtheorem{corollary}{\bf Corollary}

\newtheorem{lemma}{\bf Lemma}

\newtheorem{definition}{\bf Definition}

\renewcommand{\thetheorem}{\thepar.\arabic{nom}}

\renewcommand{\thelemma}{\thepar.\arabic{nom}}

\renewcommand{\thecorollary}{\thepar.\arabic{nom}}

\renewcommand{\thedefinition}{\thepar.\arabic{nom}}

\baselineskip0.6cm

\renewcommand{\Box}{\sharp}




{\it For Anatolii Vladimirovich YAKOVLEV 

to 60th birthday}

\

\section{Introduction}
\setcounter{par}{1}
The notion of multilinear quantum Lie operation
appears naturally in connection with a different attempts to generalize
the Lie algebras.
There are a number of reasons why the generalizations are 
necessary. First of all this is the 
demands for a "quantum algebra" which was formed in the papers by 
Ju.I. Manin, V.G. Drinfeld, S.L. Woronowicz, G. Lusztig, L.D. Faddeev,
and many others. These demands are defined by 
a desire to keep  the intuition of 
the quantum mechanics differential calculus that is based
on the fundamental concepts of the Lie groups and Lie algebras theory. 

Normally the quest for definition of bilinear brackets on the module of 
differential 1-forms that replace the Lie operation
leads to restrictions like multiplicative skew-symmetry of
the quantization parameters
\cite{Tcv}, involutivity of braidings, or 
bicovariancy of the differential calculus
\cite{Vor1}, \cite{Vor2}.
At the same time lots of quantizations, for example the Drinfeld--Jimbo
one, are defined by multiplicative non skew-symmetric parameters,
and they define not bicovariant (but one-sided covariant)
calculus. 
By these, and of course, by many other reasons, 
the 
attention of researches has been extended on operations 
that replace the Lie brackets, but depend on greater 
number of variables. Such are, for example, 
$n$-Lie algebras introduced by V.T. Filippov
\cite{Fil} and then independently appeared under the name 
Nambu--Lie algebras in theoretical
researches on generalizations of Nambu mechanics 
\cite{Dab1}, \cite{Dab2}. Trilinear operation have been
considered in the original Y. Namby paper
\cite{Dab3}, and also in a number of others papers
on generalization of quantum mechanics
(see, for example, in \cite{Yam1} a trilinear oscillator,
or in \cite{Yam2}, \cite{Yam3} multilinear commutator).
Of course one should keep in mind that both multivariable and partial
operations are the subject of investigation in the theory of
algebraic systems located at the interfaces between
algebra and mathematical logic, the algorithms theory and 
the computer calculations \cite{Mal1}, \cite{Mal2}.

Another group of problems (which is close to the author interests)
that needs the generalization of Lie algebras is connected with
researches of automorphisms and skew derivations of 
noncommutative algebras. 
A noncommutative version of the fundamental Dedekind
algebraic independence lemma says
that all algebraic dependencies in automorphisms and ordinary derivations
are defined by their algebraic structure 
(that is  a structure of a group acting on a Lie algebra)
and by operators with "inner" action
(see \cite{Khar1,IN}, Chapter 2). 

The quest for extension this, very good working, result into
the field of skew derivations
inevitably leeds to a question what algebraic structure corresponds to the
skew derivation operators? This in tern forces to consider $n$-ary
multilinear (partial) operations  on the Yetter--Drinfeld modules, 
that are not reduced to the bilinear ones
(see \cite{Khar1}, section 6.5 or   \cite{IN}, section 6.14).

In fact, the systematic investigation of quantum Lie operations 
related to the Freiderichs criteria
 \cite{Fr} have been started in the papers of B. Pareigis
\cite{Par1}, \cite{Par2}, \cite{Par3}, where he  has found 
a special series of the $n$-ary operations. Then it was continued 
in  the author papers
\cite{Khar2}, \cite{Khar3}, \cite{Khar4},
where a criteria for existence of nontrivial
quantum Lie operations is found and a version 
of the Poincar\`e--Birkhoff--Witt theorem for a class of
algebras defined by the quantum Lie operations and 
having realization inside of Hopf algebras is proposed. 

In the present paper by means of 
\cite{Khar3} we show that under the existence condition
the dimension of the space of all $n$-linear quantum Lie operations
is included between 
 $(n-2)!$ and $(n-1)!.$ The  lower bound is  achieved
if the intersection of all comforming 
(that is satisfying the existence condition)
subsets of a given set of "quantum" variables is nonempty,
while the upper bound does if the quantification matrix 
is multiplicative skew symmetric or, equivalently, all
subsets are conforming. In the latter case, as well as 
in the case of ordinary Lie algebras, all
$n$-linear operations are superpositions of the only 
bilinear quantum Lie operation, the colored super bracket.

It is interesting to note that even for   $n=4$
not all values of the mentioned above interval
are achieved for different values of quantization parameters.
Indeed, the dimension of the quadrilinear operations space may
have values 2, 3, 4, 6, while it is never equal to 5
(see the proof of the second part 
of Theorem 8.4 \cite{Khar2}).

In the last section we show that almost always 
the quantum Lie operations space is generated by
symmetric ones, provided that the ideal of the quantization parameters
is invariant with respect to the permutation group action
(only in this case the notion of a {\it symmetric operation} make sence).
 We show also all possible 
exceptions.
The space of "general" $n$-linear quantum Lie operations 
is not an exception, that is this space is always generated by 
symmetric "general" quantum Lie operations.

\section{Preliminaries}
\stepcounter{par}
Recall that a {\it quantum variable} is a variable 
$x,$ with which an element $g_x$ 
of a fixed Abelian group 
$G$ and a character $\chi :G\rightarrow {\bf k}^*$ 
are associated, where ${\bf k}$ is a ground field.
A set of quantum variables is said to be 
{\it conforming}, if
$\prod _{1\leq i\neq j\leq n}p_{ij}=1,$ where
$||p_{ij}||$ is a matrix of quantization parameters: 
$p_{ij}=\chi ^{x_i}(g_{x_j}).$
A {\it quantum operation} in quantum variables
$x_1, \ldots , x_n$ (see \cite{Khar2}) 
is a noncommutative polynomial in these variables
that has skew primitive values in every
Hopf algebra $H,$ provided that $H$ contains the group Hopf algebra
${\bf k}[G]$ 
and every variable $x_i$ has a skew primitive semi-invariant value
$x_i=a_i$:
\begin{equation}
\Delta (a_i)=a_i\otimes 1+g_{x_i}\otimes a_i;
\ \ g^{-1}a_ig=\chi ^{x_i}(g)a_i
\label{AI}
\end{equation}
for all $g\in G.$ A nonzero $n$-linear quantum operation exists
if and only if the set
$x_1,\ldots ,x_n$ is conforming, \cite{Khar3}. All the operations 
have a commutator representation
\begin{equation}
{\bf f}(x_1,\ldots ,x_n)=
\sum _{\nu \in S^1_n}\beta _{\nu }
[\ldots [[x_1,x_{\nu (2)}],x_{\nu (3)}],\ldots ,x_{\nu (n)}],
\label{pred}
\end{equation}
where  $S^1_n$ is the permutation group of $2, \ldots, n,$
while the bracket is a skew commutator 
$[u,v]=uv-p(u,v)vu$ 
with the bimultiplicative coefficient
$p(u,v).$ This coefficient is defined on the set 
of homogeneous polynomials by means of the
quantization matrix: $p(x_i,x_j)=p_{ij},$
see \cite{Khar2}. The skew commutator
$[x_i,x_j]$ is a quantum operation if and only if 
$p_{ij}p_{ji}=1.$
In particular if the quantization matrix is multiplicative
skew symmetric, then all multilinear polynomials of the form
(\ref{pred}) are quantum operations. Therefore in this particular case
the dimension of the space of operations equals $(n-1)!.$

In general, by Theorem 3 \cite{Khar3},
the linear combination (\ref{pred}) is a quantum operation if
and only if the coefficients $\beta _{\nu }$ satisfy the following
system of equations
\begin{equation}
\sum _{\nu \in N^1(s)}\beta _{\nu \mu}t^{\mu }_{\nu,s}=0,\ \ \ 
\mu \in S_n^1, 1<s<n. 
\label{ur}
\end{equation}
Here $N^1(s)$ is the set of all $s$-shuffle
from $S^1_n,$ that is the set of all permutations
$\nu $ with
\begin{equation}
\nu ^{-1}(2)<\nu ^{-1}(3)<\ldots <\nu ^{-1}(s);\ \ 
\nu ^{-1}(s+1)<\nu ^{-1}(s+2)<\ldots <\nu ^{-1}(n).
\label{N}
\end{equation}
The coefficients  $t^{\mu }_{\nu ,s}$ 
are particular polynomials in $p_{ij}, p_{ij}^{-1},$
precisely, $t^{\mu }_{\nu ,s}=\varphi (T_{\nu ,s}^{\mu }),$
see formula (\ref{V2}) below. The above system of equations
can be rewritten in the form of relations in a crossed product.

\section{Crossed product}
\stepcounter{par}
Consider a free Abelian group
${\cal F}_n,$ 
freely generated by symbols 
 $P_{ij},$\  $1\leq i\neq j\leq n.$
Denote by ${\cal P}_n$ a group algebra of this group 
over the minimal subfield  ${\bf F}$ of the field ${\bf k}.$
Clearly, ${\cal P}_n$ has a field of fractions
${\cal Q}_n$  that is isomorphic to the field 
of rational functions $\hbox{\bf F}(P_{ij}).$
The elements $P_{ij}\ ,$\ $1\leq i\neq j \leq m$ 
generate a subalgebra ${\cal P}_m$ of ${\cal P}_n.$
The action of the symmetric group $S_n$
is correctly defined on the ring ${\cal P}_n$ 
and on the field ${\cal Q}_n$ 
by $P_{ij}^{\pi }=P_{\pi (i)\pi (j)}.$ 
Thus we can define a crossed product
${\cal Q}_n*S_n$ (with a trivial factor-system).
This crossed product is isomorphic to the algebra of 
all $n!\times n!$ matrices over the Galois field 
${\cal Q}_n^{S_n},$ and it contains the skew group algebra
${\cal P}_n*S_n.$ Recall that in the trivial crossed product
the permutations commute with coefficients according to the formula
$A\pi =\pi A^{\pi }$ (see
\cite{Skr}, \cite{Ja}, \cite{Her} or other textbooks in ring theory). 

If the parameters $p_{ij}$ are defined by the quantum variables
$x_1,$ $\ldots ,x_n,$ then there exists a uniquely defined homomorphism  
\begin{equation}
\varphi : {\cal P}_n\rightarrow \hbox{\bf k},\ \ \varphi (P_{ij})=p_{ij}.
\label{gom}
\end{equation}

 If $A\in {\cal F}_n$ then by $\overline{A}$ we denote a word appearing
from  $A$ by replacing all letters $P_{ij}$ with $P_{ji}.$
We call the words $A$ and  $\overline{A}$ {\it conjugated}.
We define  
\begin{equation}
\{ A\} =A-\overline{A}^{-1}.
\label{-}
\end{equation}
For two arbitrary indices $m,\, k$ we denote by  $[m;k]$
a monotonous cycle starting with $m$ up to $k:$
\begin{equation}
[m;k]
{\buildrel \rm def \over =}
\left\{ \matrix{(m,m+1,\ldots ,k),&\ \ \hbox{if}\ m\leq k \cr 
(m,m-1,\ldots ,k),&\ \ \hbox{if}\ m\geq k.}\right.
\label{tc}
\end{equation}
Clearly $[m;k]^{-1}=[k;m]$ in these denotations.
It is easy to see that a permutation $\nu $ belongs $N^1(s)$
if and only if 
\begin{equation}
\nu =[2;k_2][3;k_3]\cdots [s;k_s]
\label{kas}
\end{equation}
for a sequence of indices
$1<k_2<k_3<\ldots k_s\leq n$ (see Lemma 1 in \cite{Khar3}).

Let us fix the following denotations for particular elements
of the skew group algebra
\begin{equation}
V_s=\sum _{\nu \in N^1(s)}\nu T_{\nu ,s}=\sum _{1<k_2<k_3<\ldots <k_s\leq n} 
[2;k_2][3;k_3][4;k_4]\cdots [s;k_s]T_{k_2:k_3:\ldots :k_s},
\label{V1}
\end{equation}
where 
\begin{equation}
T_{\nu ,s}=T_{k_2:k_3:\ldots :k_s}=
\{ \prod _{m=2}^{m=s} 
(P_{1\ m}\cdots P_{m-1\ m}\cdot P_{s+1\ m}\cdots P_{s-m+k_m\ m}) \} .
\label{V2}
\end{equation}
The braces are defined in (\ref{-}).
\addtocounter{nom}{1}
\begin{theorem}$\!\!\!.$ If $\prod _{i\neq j}p_{ij}=1,$
then $(\ref{pred})$ is a quantum operation if and only if
\begin{equation}
\sum _{\nu }\beta _{\nu }\nu ^{-1} \cdot V_s\in \mbox{ker}(\varphi )S_n^1
\label{be}
\end{equation}
for all  $s,\ 2\leq s<n,$
where $V_s$ are defined in $(\ref{V1}).$
\label{3.3}
\end{theorem}

PROOF. See Theorem 4 \cite{Khar3}. \hfill$\Box$

Let $\Sigma $ be an arbitrary multiplicative 
subset of  ${\cal P}_n$ that does not intersect ${\rm ker}(\varphi ).$ 
Consider a localization (ring of quotients) 
${\cal P}_n\Sigma ^{-1}.$ The homomorphism
$\varphi $ has a unique extension up to 
a homomorphism of 
${\cal P}_n\Sigma ^{-1}$ into the field ${\bf F}(p_{ij})$ via
\begin{equation}
\varphi (P\sigma ^{-1})=\varphi (P)\varphi (\sigma )^{-1},\ \ \ 
P\in {\cal P}_n.
\label{Si}
\end{equation}
\addtocounter{nom}{1} 
\begin{theorem}$\!\!.$ If $\prod _{i\neq j}p_{ij}=1$
and  $\Sigma $ is an arbitrary multiplicative subset of ${\cal P}_n$
that do not intersect ${\rm ker}(\varphi ),$
then an element  $B=\sum_{\nu \in S_n^1}B_{\nu }\nu ^{-1}$
with coefficients in  ${\cal P}_n\Sigma ^{-1}$
defines a quantum operation
$\sum _{\nu \in S^1_n}\varphi (B_{\nu })
[\ldots [[x_1,x_{\nu (2)}],x_{\nu (3)}],\ldots ,x_{\nu (n)}]$
if and only if 
\begin{equation}
B \cdot V_s\in \mbox{ker}(\varphi )\Sigma ^{-1}S_n^1
\label{beS}
\end{equation}
for all  $s,\ 2\leq s<n.$
\label{3.3S}
\end{theorem}
 
PROOF. It is enough to multiply $B$ from the left by a common denominator
$\sigma \in \Sigma $ of all coefficients and to apply Theorem
\ref{3.3}. \hfill$\Box$

Recall that a {\it conforming ideal} is an ideal $I$
of the algebra ${\cal P}_n,$ generated by all elements of the form 
 $\{ W\} ,$ where $W$ is an arbitrary semigroup
word in $P_{ij}$ of length $C_n^2$ that has neither 
double nor conjugated letters.
It is easy to see that the variables $x_1,\ldots ,x_n$
are conforming if and only if the ideal  ${\rm ker} (\varphi )$
contains the conforming ideal. It is very important to note that the 
conforming ideal $I$ is invariant with respect to the action
of $S_n$ (unlike the ideal ${\rm ker} (\varphi )$ in general). Therefore
the two-sided ideal of ${\cal P}_n*S_n$ generated by $I$ coincides
the right ideal $IS_n.$ 

Consider a field of rational functions ${\bf K}$ over
${\bf F}$ in $C_n^2-1$ variables
$\overline{P}_{ij},$\ $1\leq i\neq j\leq n,$\
$(i,j)\neq (1,n).$ 
Denote
\begin{equation}
\overline{P}_{1n}=(\prod_{(k,l)\neq (1,n)}\overline{P}_{kl})^{-1}. 
\label{Ksi}
\end{equation}
In this case the map $\xi : P_{ij} \rightarrow \overline{P}_{ij}$
defines an embedding of ${\cal P}_n/I$ in ${\bf K}.$

Consider a new set of quantum variables
$X_1, X_2, \ldots ,X_n$
with which free generators $G_1,\ldots ,G_n$
of a  free Abelian group  ${\bf G}$ are associated, while
the characters over $\bf K$ are defined by
 $\chi ^{X_i}(G_j)=\overline{P}_{ij}.$
\addtocounter{nom}{1}
\begin{definition}$\!\!.$ \rm 
The system (\ref{ur}) with $p_{ij}:=\overline{P}_{ij}$
is said to be the {\it general basic} system. Its solutions 
define {\it general} quantum operations in $X_1,X_2,\ldots ,X_n$
with coefficients in {\bf K}. 
\label{OBx}
\end{definition}

Denote by $\Sigma _0$ a set of all elements
$f\in {\cal P}_n$ that do not belong to $I.$ 
This is multiplicative set. By Theorem
 \ref{3.3S} an element  $B$ corresponds to a general
quantum operation if and only if
$B\cdot V_s\in I\Sigma _0^{-1}*S_n^1.$ If it is possible to 
extend the homomorphism $\varphi $ to the coefficients 
of  $B,$ then the general quantum operation defines a quantum operation
over the ground field.
\addtocounter{nom}{1}
\begin{lemma}$\!\!.$
Let $B$ corresponds to a general quantum operation.
If $\prod _{1\leq i\neq j\leq n}p_{ij}=1$ and the coefficients 
$B_{\nu }$ of $B$ belong to ${\cal P}_n\Sigma ^{-1}$
for a set $\Sigma $ that does not intersect 
${\rm ker}(\varphi ),$ than $B$ defines a quantum operation
in $x_1,\ldots ,x_n$ over ${\bf F}(p_{ij}).$
\label{OBch}
\end{lemma}

PROOF. By theorem \ref{3.3S} it is enough to show that
 $B\cdot V_s\in I\Sigma ^{-1}S_n^1.$
The left hand side of the above inclusion belongs
$I\Sigma _0^{-1}S_n^1\cap {\cal P}_n\Sigma ^{-1}S_n^1.$
This intersection equals the right hand side. Indeed,
if $i\sigma _0^{-1}=P\sigma^{-1},$ then $i\sigma =P\sigma _0.$
Since $I$ is simple ideal and $\sigma _0\notin I,$ we get $P\in I.$
 \hfill$\Box$

\section{Decreasing modules}
\stepcounter{par}
Recall that by $S_n^{k,l,\ldots }$
we denote a group of permutations $\nu \in S_n$ with
$\nu (k)=k, \nu (l)=l,\ldots .$ A right module over
${\cal P}_n*S_n^{1,n}$ generated by elements
$A_2,$\ $\ldots $\ $A_n$
is called a {\it decreasing module}, provided that the generators are 
connected by the following relation
\begin{equation}
A_kV_{(k)}+A_{k-1}D_{(k)}=0,\ \ \ 3\leq k\leq n,
\label{mod}
\end{equation}
where  $V_{(k)}$ and $D_{(k)}$ are defined by
\begin{equation}
V_{(k)}=\sum _{l=2}^{k-1}[2;l]T_{[2;l]}.
\label{vv}
\end{equation}
\begin{equation}
D_{(k)}=\tau _kT_{[2;n]}^{(2,n)}+\sum _{l=2}^{n-1}[2;l]T^{(k,n)}_{[2;l]}.
\label{d1}
\end{equation}
Here $\tau _{k},$ is set by
\begin{equation}
\tau _{k}=\left\{
\matrix{
[2;k-1][k;n-1], \hfill  &\hbox{if}\ 2<k< n,\hfill \cr
[2;n-1], \hfill  &\hbox{if}\ k=n,2 \hfill } \right. 
\label{ta}
\end{equation}
(see Definition  8 and formulae (35), (39), (40) in \cite{Khar3}).

For an arbitrary 
\begin{equation}
B={\cal A}_2(2,n)+{\cal A}_3(3,n)+\ldots +
{\cal A}_{n-1}(n-1,n)+{\cal A}_n\hbox{id},
\label{B1}
\end{equation}
where ${\cal A}_i\in {\cal P}_n*S_n^{1,n},$ 
the following formula is correct.
\begin{equation}
B\cdot V_2=D_2(2,n)+
\sum _{k=3}^n({\cal A}_kV_{(k)}+{\cal A}_{k-1}D_{(k)})(k,n),
\label{OGO1}
\end{equation}
where 
\begin{equation}
D_{2}{\buildrel \rm def \over =}
{\cal A}_n[2;n-1]T^{(2,n)}_{[2;n]}+
\sum _{l=2}^{n-1}{\cal A}_l[2;l]T^{(2,n)}_{[2;l]}.
\label{k2}
\end{equation}
Therefore the equality $B\cdot V_2\in {\rm ker}(\varphi )S_n^1$
means that both $D_2\equiv 0$ in the quotient module
${\cal P}_n*S_n^1/{\rm ker}(\varphi )S_n^1,$
and 
${\cal A}_2, {\cal A}_3,\ldots ,{\cal A}_n$
generate a decreasing module in this quotient module 
(see details in \cite{Khar3}, section 5).

Note that the element $V_{(k+1)}$ equals the element $V_2$
that is defined by $x_1, x_2,\ldots ,x_k.$ 
Therefore the following statement shows in particular
that $D_2$ equals zero in the quotient module 
${\cal P}_n*S_n^{1,n}/I(\varphi )S_n^{1,n},$
provided that the elements $V_2$ for the proper 
subsets are invertible.
\addtocounter{nom}{1}
\begin{theorem}$\!\!${\rm (Theorem 6 \cite{Khar3}).}
Let $n\geq 3$ and
\begin{equation}
X=(\prod _{n>i>j>1}P_{ij})(\prod _{i=1}^{n-1}
P_{in})(\prod _{j=2}^{n-1}P_{1j}).
\label{X10}
\end{equation}
In every decreasing module the following relation is correct 
\begin{equation}
D_{2}\prod _{k=0}^{n-3}[n-1;2]V_{(n-k)}=
A_2[n-1;2]\prod _{k=0}^{n-3}[2;2+k]V_{(n-k)}^{(1,n)[3;n-1]^k}\cdot 
\{ X\} ^{[2;n-1]}(-1)^n,
\label{ww}
\end{equation}
where $D_{2}$ is defined by $(\ref{k2})$ by replacing ${\cal A}$
with $A.$
\label{5.3}
\end{theorem}

\section{ Second subcomponents}
\stepcounter{par}
Suppose now that the intersection of all conforming subsets
is not empty.
Without loss of generality we may suppose that $x_n$
belongs to all conforming subsets. For every subset
$Y\subseteq \{ x_1,x_2,\ldots ,x_n\} $
we denote by $W(Y)$ the set of all elements
of the form  $\{ W\} ,$ where $W$ is a semigroup word
in $P_{ij},$\ $x_i,x_j\in Y$ of the length 
$C_{|Y|}^2$ that does not contain neither
double nor conjugated letters. The set  $Y$
is conforming if and only if
$W(Y)\subseteq {\rm ker}(\varphi ).$
Denote by $\Sigma _m,$\ $m\leq n$ a multiplicative set generated by
$W(Y),$ where $Y$ runs through collection of all
subsets that do not contain no one of the letters
$x_m,$ $x_{m+1},$ $\ldots ,$ $x_n.$ 
The set  $\Sigma _m$ is invariant with respect to 
the action of the subgroup $S_{m-1}^1.$ This set does not intersect 
 ${\rm ker(\varphi )},$ since
${\rm ker(\varphi )}$ is a simple ideal and $x_n$
belongs to all conforming subsets.
We may define both a localization
 ${\cal P}_n\Sigma _m^{-1}$
and a skew group ring
${\cal P}_n\Sigma _m^{-1}*S_{m-1}^1,$
that are contained in the crossed product.
\addtocounter{nom}{1}
\begin{lemma}$\!\!.$ If $x_n$ belongs to all conforming subsets,
then the elements  $V_{(m)},$\ \ $3\leq m\leq n,$ defined by
$(\ref{vv})$ are invertible in ${\cal P}_n\Sigma _m^{-1}*S_{m-1}^1.$
\label{Obr}
\end{lemma}
PROOF. Let us use the induction by $m.$ By definition the set
$\Sigma _3$ contains $\{ P_{12}\} .$ Therefore 
$V_{(3)}^{-1}=$ $\{ P_{12}\} ^{-1} \cdot id \in $ 
${\cal P}_n\Sigma _3^{-1}*S_2^1.$

Let each of 
$V_{(k)},$\ \ $3\leq k\leq m$ is invertible in
${\cal P}_n{\Sigma }_k^{-1}*S_{k-1}^1.$ This is invertible in
${\cal P}_n{\Sigma }_{m+1}^{-1}*S_m^1$as well.
Consider the set of quantum variables $x_1,\ldots ,x_m.$
The element $V_2$ for this set coincides $V_{(m+1)}.$
Let us replace denotations  $n$ by $N$ and $m$ by $n.$
It is sufficient to show that there exists an element 
$B\in {\cal P}_N{\Sigma }_{n+1}^{-1}*S_n^1$
such that $B\cdot V_2\in S_n^1$
(now $V_2=\sum _{l=2}^n[2;l]T_{[2;l]}$ and $V_{(3)},\ldots ,V_{(n)}$
are invertible in ${\cal P}_N{\Sigma }_{n+1}^{-1}*S_n^1$).
Define by induction a sequence
$A_2,\ldots ,A_n\in {\cal P}_N{\Sigma }_{n+1}^{-1}*S_n^1.$ Let
$$
A_2=(-1)^n\prod _{k=0}^{n-3}[n-1;2]V_{(n-k)}\times
$$
\begin{equation}
([n-1;2]\prod _{k=0}^{n-3}[2;2+k]V_{(n-k)}^{(1,n)[3;n-1]^k}\cdot 
\{ X\} ^{[2;n-1]})^{-1},
\label{A2}
\end{equation}
where  $X$ is defined by (\ref{X10}).
All of the factors in the parentheses are invertible in
 ${\cal P}_N{\Sigma }_{n+1}^{-1}*S_n^1$:
the  $V$'s are by the inductive suppositions,
and  $\{ X\} $ is due to this element belongs $\Sigma _{n+1}.$ 
Furthermore, let
\begin{equation}
A_k=-A_{k-1}D_{(k)}V_{(k)}^{-1},\ \ \ 3\leq k\leq n,
\label{AK}
\end{equation}
where  $D_{(k)}$ are defined by (\ref{d1}).

By right multiplication of the above equality by  
$V_{(k)},$ we get that  $A_2,\ldots ,A_{n}$ 
generate a right decreasing submodule over
 ${\cal P}_n*S_n^{1,n}.$
Therefore we may use Theorem \ref{5.3}. By 
(\ref{ww}) and (\ref{A2}) we get
$$
D_2\prod _{k=0}^{n-3}[n-1;2]V_{(n-k)}=\prod _{k=0}^{n-3}[n-1;2]V_{(n-k)}.
$$
Thus the element $D_2$ defined by (\ref{k2}) with $A$
  in place of  ${\cal A}$ equals the identity permutation.
Consider an element 
$B$  defined by  (\ref{B1}) with  $A$ in place of ${\cal A}$:
$$
B=A_2(2,n)+A_3(3,n)+\ldots +A_n\ id.
$$
Formula  (\ref{OGO1}) and definition (\ref{AK})
show that 
$$
B\cdot V_2=D_2(2,n)+\sum _{k=3}^n(A_kV_{(k)}+A_{k-1}D_{(k)})(k,n)=(2,n).
$$ 
The lemma is proved \hfill$\Box$

\section{Interval of dimensions}
\stepcounter{par}
\addtocounter{nom}{1}
\begin{theorem}$\!\!.$
If  $x_1,\ldots ,x_n$ is a conforming set of different quantum variables, 
 then the dimension of the space of all multilinear 
quantum Lie operations
in this set is  greater then or equal to $(n-2)!.$
If there exists a variable that belongs to each conforming subset then
this dimension equals  $(n-2)!.$
\label{OSNO}
\end{theorem}

PROOF. In order to prove the first part of the theorem it is enough to
show that the rank of the basic system is less then or equal to
$(n-1)!-(n-2)!=$ $(n-2)!(n-2).$ This condition is equivalent 
to all minors of the order greater then or equal to
 $(n-2)!(n-2)$ being zero. Since the minors are integer functions in the 
matrix coefficients, it is enough to show that
this condition is valid for the general basic system
(see Definition \ref{OBx}). The parameters $\overline{P}_{ij}$
of this system are connected by the only relation. Therefore
the set 
$X_1,\ldots ,X_n$ has no proper conforming subsets.
By Theorem 10 \cite{Khar3} we get the required condition. \hfill$\Box$

The second statement follows the 
invertibility of the second subcomponents.
\addtocounter{nom}{1}
\begin{lemma}$\!\!.$
If   $V_{(3)},$ $V_{(4)},$ $\ldots ,$  $V_{(n)}$
are invertible in  ${\cal P}_n\Sigma ^{-1}*S_n^{1,n},$ where
$\Sigma ={\cal P}_n\setminus {\rm ker}(\varphi ),$
then the dimension of the space of all quantum Lie operations 
is less than or equal to $(n-2)!.$
\label{obrat}
\end{lemma}
Let the element $B\sum B_{\nu }\nu ^{-1}$ corresponds to a 
quantum Lie operation
(\ref{pred}) $\varphi (B_{\nu })=\beta _{\nu }$. Then 
$B\cdot V_2\equiv 0.$ By the decomposition of 
$B$ in cosets  (\ref{B1}), we see that the elements
${\cal A}_i$ generate a right decreasing submodule in
${\cal P}_n*S_n^{1,n}/{\rm ker}(\varphi )S_n^{1,n}.$ 
Therefore 
\begin{equation}
B\equiv 
{\cal A}_2 \sum _{k=2}^n(-1)^k(\prod _{i=3}^k D_{(i)}V_{(i)}^{-1})(k,n).
\label{Bmu}
\end{equation}
Thus,  the composition 
${\bf f}\rightarrow B\rightarrow {\cal A}_2
\stackrel{\varphi }\rightarrow {\bf F}(p_{ij})S_n^{1,n},$
where $B=\sum \beta _{\nu } \nu ^{-1}$ (see (\ref{pred})),
is a linear transformation with zero kernel of the   
space of multilinear quantum Lie operations  
into the space  of left linear combinations
$\sum _{\mu \in S_n^{1,n}} \alpha _{\mu }\mu $ over 
${\bf F}(p_{ij}).$ Therefore the dimension is less than or equal to
the number of elements in $S_n^{1,n}.$ The theorem is proved. \hfill$\Box$

\section{Symmetric quantum Lie operations \\
and symmetric sets of variables}
\stepcounter{par}
An operation
 $\lbrack\!\lbrack x_1,\ldots ,x_n\rbrack\!\rbrack $
is called {\it symmetric} (or {\it skew symmetric})
if for every permutation $\pi \in S_n$
the following equality is valid
\begin{equation}
\lbrack\!\lbrack x_{\pi (1)},\ldots ,x_{\pi (n)}\rbrack\!\rbrack
=\alpha _{\pi }\lbrack\!\lbrack x_1,\ldots ,x_n\rbrack\!\rbrack ,
\label{sym}
\end{equation}
where $\alpha _{\pi }\in {\bf k}.$ 
In the case of quantum operations, 
as well as in the case of arbitrary partial operations,
we have to explain what does it mean the left hand side of
 (\ref{sym}). Strictly speaking the left hand side is defined
only if  $x_{\pi (i)}$ has the same parameters
$\chi , g$ as  $x_i$ does. By definition only in this case
the substitution $x_i:=x_{\pi (i)}$ is admissible. 
In other word all the parameters
$p_{ij}$ should be equal to the same number $q.$ 
This is very rigid condition. It  excludes 
both the colored superbrackets and the Paregis operations with 
the above defined general ones.

However, we may suppose that 
$\lbrack\!\lbrack x_1,\ldots ,x_n\rbrack\!\rbrack $
is a polynomial whose coefficients depend on the quantization parameters,
 $\chi ^{x_i}, g_{x_i},$ that is there are shown
distinguished entries of $p_{ij}$ in the coefficients.
Then a substitution $x_i\rightarrow y$ means not only 
the substitution of the variable but also one of the parameters
$g_i\rightarrow g_y,$ $\chi ^{x_i}\rightarrow \chi ^y.$
In particular, the permutation of parameters in
(\ref{sym}) means the application of this permutation
to all indices: 
$p_{ij}\rightarrow p_{\pi (i)\, \pi (j)}.$

This interpretation of the equality (\ref{sym}) will be unique and
uncontradictable only if  the application of the permutation
is independent of the way how the coefficients
$\lbrack\!\lbrack x_1,\ldots ,x_n\rbrack\!\rbrack $
are represented as rational functions or polynomials in
$p_{ij}, p_{ij}^{-1}.$ 
The action of permutations is independent of the above 
representation if (and only if) 
ker$(\varphi )$ is invariant ideal with respect to the action of
$S_n.$
\addtocounter{nom}{1}
\begin{definition}$\!\!.$ \rm
A collection of quantum variables 
 $x_1,\ldots ,x_n$ is said to be 
{\it symmetric}, if ${\rm ker}(\varphi )$
is an invariant ideal with respect to  $S_n$  or, equivalently,
the formula $p_{ij}^{\pi }=p_{\pi (i)\pi (j)}$
correctly defines the action of
$S_n$ on the ring ${\bf F}[p_{ij}].$
\label{symp}
\end{definition}

Note that the symmetricity of the collection has nothing to do with the 
symmetricity of the quantization matrix, $||p_{ij}||,$ 
while it means the symmetrisity of
relations between the parameters.

Thus, in order to give a sense to the term 
"symmetric operation" we should first suppose that 
the coefficients of the operation belong to the field
${\bf F}(p_{ij}),$ which does not affect the generality,
(see \cite{Khar3}, p. 194); and then 
we should consider only symmetric sets of quantum variables.
This is a restrictable condition, yet. 
Nevertheless this condition
excludes no one of the above examples.
 Moreover, the existence of the symmetric set 
$X_i$ with cover parameters
$\overline{P}_{ij}\rightarrow p_{ij}$
is a key argument of the 
proof of both the existence theorem and its corollaries.
Therefore the symmetric collections of variables
are of special interest.

Consider a symmetric polynomial over
${\bf F}(p_{ij}):$
\begin{equation}
{\bf f}(x_1,\ldots ,x_n) =
\sum \gamma _{\mu }x_{\mu (1)}\cdots x_{\mu (n)}.
\label{symm}
\end{equation}
Without loss of generality (if necessary by applying a permutation)
 we may suppose that the monomial
$x_1x_2\cdots x_n$
has a coefficient 1.
Let us compare coefficients at
$x_{\pi (1)}x_{\pi (2)}\cdots x_{\pi (n)}$
in the both sides of (\ref{sym}). Taking into account
 (\ref{symm}), we immediately get
 $\alpha _{\pi }=\gamma _{\pi ^{-1}}^{\pi }.$
Afterwards the equality (\ref{sym}) takes a form
$$
\gamma _{\pi ^{-1}}^{\pi }\sum _{\mu \in S_n}\gamma _{\mu }
x_{\mu (1)}\cdots x_{\mu (n)}=
(\sum _{\mu \in S_n}\gamma _{\mu }x_{\mu (1)}\cdots x_{\mu (n)})^{\pi }
$$
\begin{equation}
=\sum _{\mu \in S_n}\gamma _{\mu }^{\pi }x_{\pi (\mu (1))}\cdots 
x_{\pi (\mu (n))}
=\sum _{\nu \in S_n}\gamma _{\nu \pi ^{-1}}^{\pi }
x_{\nu (1)}\cdots x_{\nu (n)}.
\label{cotc}
\end{equation}
From here $\gamma _{\mu \pi ^{-1}}^{\pi }=
\gamma _{\pi ^{-1}}^{\pi }\gamma _{\mu }.$ 
Let us replace  $\nu =\pi ^{-1}$ and then apply $\nu $ 
to the both sides of
the above equality. We get that the normed polynomial
 (\ref{symm}) is symmetric if and only if 
\begin{equation}
 \gamma _{\mu \nu}=\gamma _{\mu }^{\nu }\gamma _{\nu },
\ \ \hbox{with}\ \
\alpha _\pi =\gamma _{\pi ^{-1}}^{\pi }=\gamma_{\pi }^{-1}.
\label{cot}
\end{equation}
To put it another way, the set of normed symmetric polynomials
can be  identified with the first cogomology group 
$H^1(S_n, {\bf F}(p_{ij})^*)$
with values in the multiplicative group of 
${\bf F}(p_{ij}).$

Now a natural question arises: is the space of multilinear 
quantum Lie operations generated by the symmetrical ones, 
provided that the variables form a symmetric set?

We start with some examples.

EXAMPLES. Let the set of variables is 
{\it absolutely symmetric}, that is 
$p_{ij}=q.$ 
In this case the group $S_n$ acts identically on the 
field
{\bf F}$(p_{ij}).$ 
Therefore there exists only two symmetric polynomials:
\begin{equation}
S(x_1,\ldots ,x_n)=\sum _{\pi \in S_n}x_{\pi (1)}\cdots x_{\pi (n)},
\label{symet}
\end{equation}
\begin{equation}
T(x_1,\ldots ,x_n)=\sum _{\pi \in S_n}(-1)^{\pi }
x_{\pi (1)}\cdots x_{\pi (n)}.
\label{cosym}
\end{equation}
On the other hand, if the existence condition, 
$q^{n(n-1)}=1,$ is valid then, according to Theorem
 \ref{OSNO}, the dimension $l$ of the multilinear 
operations space can not be less than 
 $(n-2)!.$ Thus, if $n>4,$ or if $n=4$ and 
the characteristic of the ground field
equals 2, then wittingly  the basis consisting of the 
symmetric operations does not exist.

If  $n=4$ then we may use the analysis
from   \cite{Khar2} (see the proof 
of Theorem 8.4): $l=2$ only if
$q^{12}=1,$ $q^{6}\neq 1,$ $q^4\neq 1,$ or, equivalently,
$q^6=-1,$ $q^2\neq -1.$  
If under these conditions the polynomials $S,$ $T$ 
are quantum operations, then they should be expressed trough
the main quadrilinear operation 
(see  \cite{Khar2}, formula (57)) with the coefficients equal to
ones at
$x_1x_2x_3x_4$ and $x_1x_3x_2x_4.$ That is
$$
S=\lbrack\!\lbrack x_1,x_2,x_3,x_4\rbrack\!\rbrack +
\lbrack\!\lbrack x_1,x_3,x_2,x_4\rbrack\!\rbrack ,
$$
$$
T=\lbrack\!\lbrack x_1,x_2,x_3,x_4\rbrack\!\rbrack -
\lbrack\!\lbrack x_1,x_3,x_2,x_4\rbrack\!\rbrack .
$$
The sum of these equalities shows that all coefficients
of the main quadrilinear operation at monomials corresponding
to odd permutations have to be equal to zero.
Alternatively, the explicit formula 
(56), \cite{Khar2} shows that the coefficient at
$x_1x_2x_4x_3$ equals
$$
-{{\{ p_{13}p_{23}\} }\over {\{ p_{13}p_{23}p_{43}\} }}=
-{{q^2-q^{-2}}\over{q^3-q^{-3}}}\neq 0,
$$
since $q^4\neq 1.$ Thus in this case the symmetric basis neither exists. 

If  $n=3,$ the existence condition takes up the form
$q^6=1.$ If $q\neq \pm 1,$ then there exists only one
trilinear operation up to a scalar multiplication,
and this operation is symmetric
(see Theorem  8.1 and formula (46), \cite{Khar2}).
While if $q=\pm 1$ then the operation space 
is generated by two polynomials: 
$[[x_1,x_2],x_3]$ and $[[x_1,x_3],x_2].$
If both $S, T$ are quantum operations, then 
as above we get a contradiction $S+T=2[[x_1,x_2],x_3].$
\addtocounter{nom}{1}
\begin{lemma}$\!\!.$ Let the quantization matrix of a symmetric
quadruple of quantum variables has the form 
\begin{equation}
\vert \vert p_{ij}\vert \vert =
\left(\matrix{
* & p & q & s\cr
p & * & s & q\cr
q & s & * & p\cr
s & q & p & *\cr}\right),
\label{mat}
\end{equation}
where  $p,q,s$ are pairwise different and $p^2q^2s^2=1$.

$1.$ If the characteristic of  the field ${\bf k}$
is not equal to $2$ then there do not exist 
nonzero quadrilinear symmetric quantum Lie operations at all.

$2.$ If the characteristic equals $2$ then there exist not more then 
two linearly independent quadrilinear symmetric operations.

$3.$ In both cases the dimension of the
whole space of quadrilinear quantum Lie operations equals three.
\label{raz}
\end{lemma}

PROOF. If the parameter matrix has the  form (\ref{mat}),
then the action of the group $S_4$  on the field ${\bf F}(p_{ij})$ 
is not faithful. The kernel of this action wittingly includes 
the following four elements 
\begin{equation}
\hbox{id; } a=(12)(24);\ b=(13)(24); \ c=(14)(23).
\label{elem}
\end{equation}
These elements form a normal subgroup
$H\lhd S_4$ isomorphic to $Z_2\times Z_2.$
Let 
\begin{equation}
S=\sum _{\pi \in S_4}
\gamma _{\pi }x_{\pi (1)}x_{\pi (2)}x_{\pi (3)}x_{\pi (4)}
\label{opera}
\end{equation}
be some symmetric quantum operation, $\gamma _{\rm id}=1.$
According to  (\ref{cot}) with $h=\mu =\nu \in H$ we have
$\gamma _h^2=\gamma _{h^2}=\gamma _{\rm id}=1,$
that is $\gamma _h=\pm 1\in {\bf F}.$ Moreover, all of the elements 
$\gamma _h,$ $h\neq $ id, $h\in H$ may not be equal to $-1,$ 
since, again by (\ref{cot}), the product of 
every two of them equals the third one. On the  other hand
formula (\ref{cot}) with $h\in H,$
$g\in S_4$ implies $\gamma _{g^{-1}}^g\gamma _g=1$ and
$$
\gamma _{g^{-1}hg}=\gamma _{g^{-1}}^{hg}\gamma _{hg}=
\gamma _{g^{-1}}^g\gamma _h^g\gamma _g=\gamma _h^g=\gamma _h.
$$
Therefore all of  $\gamma _h,$ $h\in H$ equal each other
and equal to 1.
 
Furthermore, the polynomial $S,$ 
as well as any other quantum Lie operation,
has a commutator representation
(\ref{pred}): 
\begin{equation}
S=\sum _{\nu \in S_4^1} 
\beta _{\nu }[[[x_1,x_{\nu (2)}],x_{\nu (3)}],x_{\nu (4)}].
\label{razl}
\end{equation}
If we compare coefficients at monomials
$x_1x_2x_3x_4$ and $x_4x_3x_2x_1,$
we get 
$1=\gamma _{\rm id}=\beta _{\rm id}$ and
$
1=\gamma _{(14)(23)}=$ $ \beta _{\rm id}(-p_{12})(-p_{13}p_{23})
(-p_{14}p_{24}p_{34})=$ $-p^2q^2s^2=-1.
$
This complete the first statement.

In both  cases the condition
$p^2q^2s^2=1$ means that all three element subsets of the given
quadruple are conforming. If some pair of
them does as well, say  $1=p_{12}p_{21}=p^2,$
then by symmetricity all others pairs are conforming too, that is 
$q^2=s^2=1.$ In this case  $p,q,s\in {\bf F}.$
Thus  $p=p^{(23)}=q=q^{(34)}=s.$ This is contradiction with 
the lemma condition. Therefore by Theorem 8.4 \cite{Khar2}
(see the second case in the proof of the second part) 
the quadrilinear quantum Lie operations space 
is generated by the following three polynomials
\begin{equation}
[W,x_4]; \ \ [W^{\sigma },x_1]; \ \ [W^{\sigma ^2}, x_2], 
\label{opp}
\end{equation}
where $\sigma =(1234)$ is the cyclic permutation,
while $W$ is the main trilinear operation in $x_1,x_2,x_3.$ 
By the definition of this operation, see \cite{Khar2} formula (45),
in the case of the characteristic 2, we get
\begin{equation}
W=(x_1x_2x_3+x_3x_2x_1)+
{{p+p^{-1}}\over{q+q^{-1}}}(x_2x_3x_1+x_1x_3x_2)
{{s+s^{-1}}\over{q+q^{-1}}}(x_3x_1x_2+x_2x_1x_3).
\label{r45}
\end{equation}
Let
\begin{equation}
S=\xi [W,x_4]+\xi _1[W^{\sigma },x_1]+\xi _2[W^{\sigma ^2}, x_2]. 
\label{str}
\end{equation}
If we compare the coefficients at the monomials
 $x_1x_2x_3x_4$
and $x_2x_1x_4x_3,$ we get
$\xi +\xi _1=\gamma _{\rm id}=1,$ 
$\xi _2=\gamma _{(12)(34)}=1.$ Therefore
\begin{equation}
S=\xi ([W,x_4]+[W^{\sigma },x_1])+
([W^{\sigma },x_1]+[W^{\sigma ^2}, x_2]). 
\label{strrr}
\end{equation}
Thus the symmetric operations generate 
not more then two-dimensional subspace. The lemma is proved.
 \hfill $\Box$  
\addtocounter{nom}{1}
\begin{theorem}$\!\!.$
If $x_1,x_2,\ldots ,x_n$ is a symmetric but not absolutely symmetric
collection of quantum variables, then the 
multilinear quantum Lie operations space
is generated by symmetric operations, with the only 
exception given in the lemma $\ref{raz}.$
\label{symett}
\end{theorem}
PROOF. Consider the skew group algebra
$M={\bf F}(p_{ij})*S_n.$ The permutation group action defines 
a structure of right $M$-module on the set of multilinear 
polynomials:
\begin{equation}
\sum \gamma _{\pi }x_{\pi (1)}\cdots x_{\pi (n)}\cdot 
\sum \beta _{\nu } \nu =
\sum _{\nu , \pi }
\beta _{\nu }\gamma _{\pi }^{\nu }x_{\nu (\pi (1))}\cdots x_{\nu (\pi (n))}.
\label{deiis}
\end{equation}
A polynomial (\ref{symm}) is symmetric if and only if 
it generates a submodule of dimension one over
${\bf F}(p_{ij})$.

Note that the quantum Lie operations space form 
a right $M$-submodule, provided that the collection 
of variables is symmetric.
Indeed, let the system (\ref{ur}) is fulfilled for the coefficients
of a polynomial ${\bf f}$ defined by (\ref{pred}).
The application of a permutation $\pi \in S_n^1$ to this system
shows that the coefficients of the polynomial ${\bf f}^{\pi }$
satisfy the same system up to rename of the variables
$x_i\rightarrow x_{\pi (i)}.$ Therefore ${\bf f}^{\pi },$
$\pi \in S_n^1$ is a quantum Lie operation. 
If we replace the roles of indices 1 with 2, we get that
${\bf f}^{\pi },$ $\pi \in S_n^2$ is a quantum Lie operation as well.
Since the subgroups  $S_n^1$ and $S_n^2$ with $n>2$ 
generate $S_n,$ all multilinear quantum Lie operations form an
$M$-submodule.
  
Suppose now that the $S_n$ action on the field
${\bf F}(p_{ij})$ is faithful. In this case  $M$ is isomorphic to
the trivial crossed product of the field ${\bf F}(p_{ij})$
with the Galois group $S_n$. Consequently 
$M$ is isomorphic to the algebra of $n!$ by $n!$
matrices over the Galois subfield
${\bf F}_1={\bf F}(p_{ij})^{S_n}.$ 
Therefore each right $M$-module equals a direct sum 
of irreducible submodules, while all irreducible submodules are
isomorphic to the $n!$-rows module over the Galois field 
${\bf F}_1.$  On the  other hand, the dimension of 
${\bf F}(p_{ij})$ over  ${\bf F}_1$ equals $n!$ too.
Since every right $M$-module is a right
space over  ${\bf F}(p_{ij}),$
all irreducible right $M$-modules 
are of dimension one over ${\bf F}(p_{ij}).$
This proves the theorem in the case of  a faithful action.

If $n>4$ or $n=3,$ while the action is not faithful,
then all even permutations act identically.
This immediately implies that the collection of variables is 
absolutely symmetric, $p_{ij}=q.$

Let $n=4.$ If the action is not faithful then all the permutations
 (\ref{elem}) act identically. This implies the parameter matrix
has the form (\ref{mat}). The existence condition for quantum 
Lie operations is $p^4q^4s^4=1.$ If  $p^2q^2s^2=1$ then
we get the lemma \ref{raz} exception. 
Therefore suppose that $p^2q^2s^2=-1\neq 1.$

If $p,q,s$ are pairwise different, then 
$S_4^1$ acts faithfully on ${\bf F}(p,q,s).$ 
Therefore $M_1={\bf F}(p,q,s)*S_4^1$ is the algebra of
$6$ by $6$ matrices over the Galois field ${\bf F}_1.$
This is central simple algebra. Thus it splits in $M$
as a tensor factor $M=M_1\otimes Z_1,$ where 
$Z_1$ is a centralizer of  $M_1$ in $M.$ Let us calculate this 
centralizer.

First of all, $Z_1$ is contained in the centralizer of 
${\bf F}(p,q,s),$ that equals the group algebra 
$A={\bf F}(p,q,s)[{\rm id},a,b,c].$ This group algebra has a 
decomposition in a direct sum of ideals
$$
A={\bf F}(p,q,s)e_1\oplus {\bf F}(p,q,s)e_2\oplus
{\bf F}(p,q,s)e_3\oplus {\bf F}(p,q,s)e_4,
$$
where  $e_1={1\over 4}({\rm id}+a+b+c),$ 
$e_2={1\over 4}({\rm id}+a-b+c),$ $e_3=e_2^{(23)},$ $e_4=e_2^{(34)}.$
The stabilizer of $e_2$ in $S_4^1$ equals a two-element
subgroup $S_4^{1,3}.$ Let ${\bf F}_2={\bf F}(p,q,s)^{S_4^{1,3}}$
be a Galois subfield of this subgroup. Then $Z_1$ equals
the centralizer of $S_4^1$ in $A.$ This consists of
sums
$$
\alpha e_1+\beta e_2+\beta ^{(23)}e_3+\beta ^{(34)}e_4,\  \ 
\alpha \in {\bf F}_1,\ \beta \in {\bf F}_2.
$$
Thus, $Z_1\simeq {\bf F}_1\oplus {\bf F}_2.$
Consequently, 
$M\simeq ({\bf F}_1)_{6\times 6}\oplus ({\bf F}_2)_{6\times 6}.$

This means that up to isomorphism there exists just two 
irreducible  right modules over $M.$ 
One of them equals the 6-rows space 
over ${\bf F}_1,$ while another one equals the 6-rows space
over ${\bf F}_2.$ The dimensions of these modules over 
${\bf F}_1$ equal to respectively  6 and 18.
Therefore, the first module is of dimension one over 
${\bf F}(p,q,s),$ while the second one is of dimension
three. By Theorem 8.4 \cite{Khar2} the quantum Lie 
operation module is of dimension two over ${\bf F}(p,q,s).$
Therefore its irreducible submodules may not be of dimension three.
Thus all of them are of dimension one. The theorem is proved.
 \hfill$\Box$
\addtocounter{nom}{1}
\begin{corollary}$\!\!.$
There exists a collection of 
$(n-2)!$ general symmetric multilinear quantum Lie operations
that generates the space of all the operations.
\label{sl1}
\end{corollary}
The same statement is valid for quantum variables considered by
Paregis in \cite{Par1} as well, that is in the case when the
 quantization parameters are related by
$p_{ij}p_{ji}=\zeta ^2,$ where $\zeta $ is a $n$th primitive 
root of 1.
\addtocounter{nom}{1}
\begin{corollary}$\!\!.$
The total number of linearly independent 
symmetric multilinear quantum Lie operations for symmetric,
but not absolutely symmetric, Pareigis quantum variables 
is greater than or equal to $(n-2)!.$
\label{sl2}
\end{corollary}

\

V.K. Kharchenko, Universidad Nacional Aut\'onoma de M\'exico, 
Campus Cuautitlan, Cuautitlan Izcalli, Estado de M\'exico, 54768, 
M\'exico and

E-mail: vlad"servidor.unam.mx

\begin{thebibliography}{93}
\bibitem{Tcv}
C. Bautista, A Poincar\`e--Birkhoff--Witt theorem for 
generalized color algebras, 
Journal of Mathematical Physics, 39, N7(1998), 3829--3843.
\bibitem{Vor1} S.L. Woronowicz, Compact matrix pseudogroups,
Communications in Math. Phys., 111(1985), 613--665.
\bibitem{Vor2} S.L. Woronowicz, Differential calculus on 
compact matrix pseudogroups,
Communications in Math. Phys., 122(1989), 125--170.
\bibitem{Fil} V.T. Filippov, $n$-Lie algebras, Sib. Math. Journal,
 26, N6(1985), 126--140.
\bibitem{Dab1} L. Takhtajan, On foundations of the generalized
Nambu mechanics, Communications in Math. Phys., 160(1994), 295--315.
\bibitem{Dab2} G. Dito, M. Flato, D. Sternheimer, L. Takhtajan,
Deformation quantization and Nambu mechanics,
Communications in Math. Phys.,183(1997), 1--22.
\bibitem{Dab3} Y. Nambu, Generalized Hamilton dynamics,
Physical Review D, 7, N8(1873), 2405--24--12.
\bibitem{Yam1} R.M. Yamaleev, Elements of cubic
quantum mechanics, Communications of JINR, Dubna, 
2-88-147, 1988, 1--11.
\bibitem{Yam2} R.M. Yamaleev, Model of 
plylinear oscillator in a space of non-integer 
quantum numbers. Communications of JINR, Dubna, 
2-88-871, 1988, 1--10.
\bibitem{Yam3} R.M. Yamaleev, Model of
polylinear Bose- and Fermi-like oscillator,
Communications of JINR, Dubna, 
2-92-66, 1992, 1--14.
\bibitem{Mal1} A.I. Mal'tcev, Algebraic systems, 
Moscow, "Nauka",  1965.
\bibitem{Mal2} A.I. Mal'tcev, Algorithms and recursive functions,
Moscow, "Nauka", 1970.
\bibitem{Khar1}
V.K. Kharchenko, Automorphisms and Derivations of associative Rings,
Kluwer Academic Publishers, Mathematics and its applications,
(Soviet Series), v.69, Dordrecht/Boston/London, 1991.
\bibitem{IN} V.K.Kharchenko, Noncommutative Galois theory,
Novosibirsk, "Nauchnaja Kniga", 1996.
\bibitem{Fr} K.O. Friedrichs, Mathematical 
aspects of the quantum theory of fields. V,
Communications in Pure and Applied Mathematics,
6(1953), 1-72.
\bibitem{Par1} B. Pareigis, On Lie algebras in braided categories,
in: Quantum Groups and Quantum Spaces, eds. R. Budzy\'nski,
W. Pusz, S. Zakrzewski, Banach Center Publications, {\bf 40} 
(1997), 139--158.
\bibitem{Par2} B. Pareigis, Skew-primitive elements of 
quantum groups and braided Lie algebras, Rings, Hopf algebras, 
and Brauer groups (Antwerp/Brussels), Lecture Notes
in Pure and Appl. Math., 197(1996), 219--238. 
\bibitem{Par3} B. Pareigis, On Lie algebras in the category 
of Yetter--Drinfeld modules, Appl. Categ. Structures, 6, N2(1998),
152--175.
\bibitem{Khar2} V.K. Kharchenko,
An algebra of skew primitive elements,
Algebra and Logic, 37, N2(1998), 101--127.
\bibitem{Khar3} V.K. Kharchenko, An existence condition for 
multilinear quantum operations,
Journal of Algebra, 217(1999), 188--228.
\bibitem{Khar4} V.K. Kharchenko, A quantum analog of the 
Poincar\`e--Birkhoff--Witt theorem, Algebra and Logic,
 38, N4(1999), 259--276, (QA/0005101).
\bibitem{Skr} A.A. Albert, 
Structure of Algebras, Amer. Math. Soc., Providence, RI(1961).
\bibitem{Ja} N. Jacobson, Structure of Rings, AMS Colloq.
Publ. vol. 37, Amer. Math. Soc., Providence, RI (1956, revised 1964).
\bibitem{Her} I. Herstein, Noncommutative Rings,
Carus Mathematical Monographs, N15, Amer. Math. Soc., 1998.
\end{thebibliography}
\end{document}